# A language theoretic analysis of combings

Sarah Rees,
Department of Mathematics, University of Newcastle,
Newcastle NE1 7RU, UK,
e-mail: Sarah.Rees@ncl.ac.uk

June 10, 1997

**Abstract**

A group is combable if it can be represented by a language of words satisfying a fellow traveller property; an automatic group has a synchronous combing which is a regular language. This paper gives a systematic analysis of the properties of groups with combings in various formal language classes, and of the closure properties of the associated classes of groups. It generalises previous work, in particular of Epstein et al. and Bridson and Gilman.

## 1 Introduction

The concept of a combing for a finitely generated group has grown out of the definition of an automatic group (as introduced in [10]); in the terminology of this paper, a group is automatic precisely when it possesses a regular synchronous combing. (Roughly speaking, a combing is an orderly set of strands through a Cayley graph of the group, which is regular if it is defined by a finite state automaton; a formal definition is given in section 2.)

The class of automatic groups has valuable properties (particularly from a computational point of view) and contains a wealth of examples. However it is disappointing for its failure to include the fundamental groups of all compact geometrisable 3-manifolds, particularly so since it was the properties of the fundamental groups of compact hyperbolic manifolds observed in work of Cannon ([9]) which motivated the definition of this class; subsequent work in [10]

The author would like to thank both the Fakultät für Mathematik of the Universität Bielefeld for its hospitality while this work was carried out, and the Deutscher Akademischer Austauschdienst for its financial support.

then showed that the fundamental groups of compact 3-manifolds modelled on 6 of Thurston's 8 geometries are in fact automatic. Further, no nilpotent group without a finite index abelian subgroup is automatic, and no soluble (non virtually abelian) example is known (at least to this author).

A result of Bridson and Gilman ([7]) suggests a way to remedy the deficiencies of automatic groups. Bridson and Gilman show that the fundamental group of any compact geometrisable 3-manifold possesses an asynchronous combing which is an indexed language (as defined in [1]); the same techniques show that the same result holds for many nilpotent and soluble groups. The class of groups with combings of this type is defined by weakening both of the two restrictions on the 'language' associated with an automatic group; the geometric 'fellow traveller condition' is relaxed from a synchronous to an asynchronous condition, and the language theoretic requirement of regularity is replaced by the requirement that the language be indexed (that is, recognised by a one-way nested stack automaton, a type of machine, defined in [2], which is a little more general than a pushdown automaton).

Most of the recent more general results on combings have concentrated on those which, as in Bridson and Gilman's case, are asynchronous; these are certainly easier to find than synchronous combings. However, from a practical point of view, synchronous combings, when they exist, are more useful. The value of the study of these is increased by a very recent, and as yet unpublished, result of Bridson and N. Brady, which constructs a synchronous combing for a group which is not biautomatic, and almost certainly not automatic.

Naturally some properties of the class of automatic groups are lost by the move to a more general definition; however, widening the definition admits more constructions within a class of combable groups, and so allows more examples. Our aim in this paper is to give a systematic analysis of the properties of various classes of synchronous and asynchronous combings defined by formal language classes, to see to what extent properties of the class of automatic groups still hold in these classes, and what further properties can be deduced. We build on results in [10], [7], [5] and [20], in particular the first two of those. Often the proofs of these papers remain valid with a small amount of adjustment in the more general context we need; details from the original proofs are sometimes repeated for greater clarity. Our goal is to prove these results in the most general context possible.

Section 2 of this article contains the definitions of various types of combings and fellow traveller properties; section 3 introduces the formal language theory which is used. Section 4 examines the solubility of the word problem, section 5 the closure of classes of combings under change of generators and various finite variations. Section 6 looks at the closure of combing classes under free and direct products, central and split extensions, and finally section 7 uses a result of Bridson, already used by Bridson and Gilman to find indexed combings for fundamental groups of *Nil* and *Sol* manifolds, to show that similar combings exist also for many other nilpotent and soluble groups (this question is examined in more detail in [13]).

## 2 Fellow traveller properties, languages and combings

Let $G$ be a group, with identity element $e$, and finite generating set $X$. Without loss of generality, we may assume that $X$ is inverse closed, that is, contains the inverse of each of its elements; we shall make this assumption thoughout this paper. We call a product of elements in $X$ a word over $X$, and denote by $X^*$ the set of all such words. Let $\Gamma = \Gamma_{G,X}$ be the Cayley graph for $G$ over $X$, with vertices corresponding to the elements of $G$, and, for each $x \in X$, a directed edge from the vertex $g$ to the vertex $gx$, labelled by $x$. Let $d_{G,X}$ measure (graph theoretical) distance between vertices of $\Gamma_{G,X}$. For words $w, v \in X^*$, we write $w = v$ if $w$ and $v$ are identical as words, $w =_G v$ if $w$ and $v$ represent the same element of $G$. We define $l(w)$ to be the length of $w$ as a string, and $l_G(w)$ to be the length of the shortest word $v$ with $v =_G w$. It is straightforward to extend $d_{G,X}$ to a differentiable metric on the 1-skeleton of $\Gamma$. Then each word $w$ can be associated with a differentiable path from $e$ labelled by $w$ such that, for $t < l(w)$, the path from $e$ to $w(t)$ has length $t$, and for $t \geq l(w)$, $w(t) = w(l(w))$.

Suppose that $v, w$ are words in $X^*$, and that $K \in \mathbf{N}$. Various *fellow traveller* properties can describe the relationship between $v$ and $w$. We say that $v$ and $w$ *synchronously* $K$-*fellow travel* if for all $t$, $d_{G,X}(v(t), w(t)) \leq K$. More generally, we say that $v$ and $w$ *asynchronously* $K$-*fellow travel* if there is a strictly increasing, differentiable function $h : \mathbf{R} \to \mathbf{R}$, mapping $[0, l(v)]$ onto $[0, l(w)]$, with the property that, for all $t > 0$, $d_{G,X}(v(t), w(h(t))) \leq K$. We say that $v$ and $w$ asynchronously $K$-fellow travel with bound $M$, if for all $t \leq l(v)$, the function $h$ above satisfies $1/M \leq h'(t) \leq M$. We shall call $h$ the *relative speed* function of $v$ and $w$, and $K$ the *fellow traveller constant*. Note that $v$ and $w$ synchronously fellow travel if and only if they asynchronously fellow travel with relative speed function $h$ defined by $h(t) = t$.

We define a *language* for $G$ over $X$ to be a set of words in $X$ which contains at least one representative for each element of $G$. A language is *bijective* if it contains exactly one representative of each group element, *prefix closed* if it contains every prefix of each of its elements. It is *geodesic* if for each $w \in L$, $l(w) = l_G(w)$, and *near geodesic* if, for some constant $\epsilon$, for all $w \in L$, $l(w) - l_G(w) < \epsilon$ (this is stronger than being quasigeodesic, see [10]).

Suppose that $L$ is a language for $G$. We call $L$ a *synchronous combing* if for some constant $K$ (the *fellow travelling constant* for $L$), the $K$-fellow traveller condition is satisfied by all pairs of words $v, w \in L$ for which $w =_G vx$ for some $x \in X \cup \{e\}$. We call $L$ a *synchronous bicombing* if the $K$-fellow traveller condition is also satisfied by pairs of words $xv, w$ with $v, w \in L$ and $w =_G xv$, for $x \in X$. We call $L$ an *asynchronous combing* if for some $K$ the asynchronous $K$-fellow traveller condition is satisfied by all pairs of words $v, w \in L$ for which $w =_G vx$ for some $x \in X \cup \{e\}$, and an *asynchronous bicombing* if the same is true for words $xv, w$ with $v, w \in L$ and $w =_G xv$, for $x \in X$. We say that an asynchronous combing (or bicombing) $L$ is *boundedly asynchronous* if for some constant $M$, relevant pairs of words asynchronously fellow travel with bound $M$. (Boundedly asynchronous combings are important for the study of asynchronously automatic groups in [10].)

When we state results, we need to be very clear about the various fellow traveller conditions which are currently under discussion. Hence we shall consider the words *combing* and *bicombing* to have a rather general meaning, when used without further qualification. We introduce the a notion of *type* of a combing for situations where we need to be more precise.

More specifically, we shall use the term *combing* for a language which is either a synchronous, asynchronous, or boundedly asynchronous combing or bicombing, and call a group *combable* if it has a combing. Similarly, we shall use the general term *bicombing* for a language which is either a synchronous, asynchronous, or boundedly asynchronous bicombing, and call a group *bicombable* if it has a combing. A *type*, $\tau$, of combing will be one of the six classes of combings (asynchronous, boundedly asynchronous or synchronous combings or bicombings) we have defined above; we shall call the set of all six classes $\mathcal{T}$, and for $\tau \in \mathcal{T}$, we shall say that a combing in the class $\tau$ *has type* $\tau$. (Clearly a combing may therefore have more than one type.) The following inclusions between these classes are fairly clear:-

$$\begin{array}{ccccc}
\text{sync. combings} & \subset & \text{bdd. async. combings} & \subset & \text{async. combings} \\
\cup & & \cup & & \cup \\
\text{sync. bicombings} & \subset & \text{bdd. async. bicombings} & \subset & \text{async. bicombings}
\end{array}$$

We call fellow traveller properties holding between pairs of words $v, w$ with $w =_G v$ or $w =_G vx$ *right* fellow traveller properties, and fellow traveller properties holding between pairs of words $v, w$ with $w =_G xv$ *left* fellow traveller properties.

We call a combing *bijective*, *prefix closed*, *geodesic* or *near geodesic* if the underlying language has those properties. Note that the literature contains various uses of the term combing; [10], for instance, requires that the language consist of quasigeodesics, and [5] and [7] require that the language be bijective. We prefer to use a more general definition. In fact, when we construct combings out of other combings, it will often be clear from the construction that properties such as bijectivity, prefix closure, geodesicity or near geodesicity possessed by the original combings would be inherited by the new ones; in general, however the construction itself is not dependent on those properties.

We have the following general results for combable groups, from [5] and [11]. (In fact both authors assumed bijectivity as part of their definition of a combing, but the condition was not used.)

**Proposition 2.1 ([5])** *Any asynchronously combable group is finitely presented.*

**Proposition 2.2 ([11])** *Any asynchronously combable group has an exponential isoperimetric inequality. A group with a prefix closed asynchronous combing has a quadratic isoperimetric inequality.*

As an immediate corollary to proposition 2.2 (see, for instance [10], theorem 2.2.5) we have

**Corollary 2.3** *Any asynchronously combable group has soluble word problem.*

Hence, since finitely presented soluble groups of derived length 3 are known with insoluble word problem ([19, 4]), the following is obvious.

**Corollary 2.4** *There exist finitely presented soluble groups of derived length 3 which are not asynchronously combable.*

A natural restriction to put on a combing is to require that it is a language in one of the families of formal languages, that is, a language recognised by some theoretical model of computation (or equivalently, defined by a formal grammar). The combing of an automatic group is required to be a regular language. Bridson and Gilman in [7] studied bijective, asynchronous combings in various formal language families, in particular the families of bounded languages, regular languages, context-free languages and indexed languages, and more generally any full abstract family of languages (of which the last three are examples, see [18] ). We aim to give a more general analysis, covering both a range of fellow traveller conditions (synchronous and even boundedly asynchronous conditions are in general harder to obtain than asynchronous ones) and a range of language families.

## 3 Families of languages

Following [18], we define a family of languages to be any non-empty set of languages over finite alphabets; hence in particular, we include under this definition the set $\mathcal{U}$ of all languages with finite alphabet. Of course, all the closure properties for a language family given as conditions for the results in the following sections are valid in $\mathcal{U}$, and hence those results hold for combings for which there is no language theoretic restriction. However, in this paper, we are chiefly interested in combings associated with families of formal languages.

Many families of formal languages are described in [18], including the four families of the Chomsky Hierarchy (regular sets, context-free languages, context-sensitive languages and recursively enumerable sets), the recursive sets, and (briefly) the indexed languages; [14] provides a more detailed mathematical treatment of context-free languages. [12] provides a very useful introduction to formal language theory from a group theoretical standpoint. The indexed languages were originally defined, and are more fully described in [1] and [2]. Bounded languages are described in [16].

A formal language may be defined by a machine which recognises it; for example, regular sets are recognised by finite state automata, context-free languages by (non-deterministic) pushdown automata, recursively enumerable and recursive sets by Turing machines and halting Turing machines respectively. Alternatively, a formal language may be defined by a set of grammatical rules. Such a description may well be more digestible; certainly this seems to be the case for indexed languages, recognised by one way nested stack automata, and generated by grammars described in [1], and also in [7].

From now on, we shall use the symbol $\mathcal{F}$ to denote a family of languages, and in particular $\mathcal{R}$ to denote the family of regular sets, $\mathcal{C}$ the context-free languages, $\mathcal{DC}$ the deterministic subfamily of $\mathcal{C}$, $\mathcal{I}$ the indexed languages, $\mathcal{R}ec$ the recursive sets and $\mathcal{R}en$ the recursively enumerable sets. We call a combing for a group $G$ an $\mathcal{F}$-combing if it is a combing in the family $\mathcal{F}$ of languages. Of course, every combing is a $\mathcal{U}$-combing.

Hopcroft and Ullman follow a general approach to the study of formal languages proposed by [15] and based on the closure properties of the various families. Recognising that some closure properties are consequences of others, they introduce the concepts of abstract and full abstract families of languages, and also of trios and full trios, for which various closure properties are satisfied. More precisely, where $P_{\cap \mathcal{R}}$ etc. are as defined in table 1, a family of languages is called a *trio* if it satisfies $P_{\cap \mathcal{R}}$, $P_{\epsilon-freeh}$ and $P_{h^{-1}}$, a *full trio* if it is a trio which also satisfies $P_h$, an *abstract family of languages* (AFL) if it is a trio which also satisfies $P_{\cup}$, $P_{\circ}$ and $P_{+}$, and a *full* AFL if it is a full trio which also satisfies $P_{\cup}$, $P_{\circ}$ and $P_{*}$. The families $\mathcal{R}$, $\mathcal{C}$, $\mathcal{I}$ and $\mathcal{R}en$ are full AFL's; $\mathcal{R}ec$ is an AFL ([18, 7]).

The properties of formal languages defined in table 1 are used in this paper (or seem to be of related interest). The results quoted are mostly taken from [18]; some results for indexed languages come from [1, 2, 17].

## 4 Recognising fellow travellers, bijectivity and the word problem

Where $G$ is a group with finite generating set $X$, for any constant $K$, the set of pairs of words $(w, v)$ in $X^*$ which synchronously $K$-fellow travel is a regular set. It is recognised by a (2-stringed) finite state automaton which is usually known as a *difference machine* $D_K$; the states of $D_K$ correspond to 'word differences' (that is, the words of length at most $K$ which connect a pair of fellow travellers), and transitions encode the equations $dx =_G x'd'$ for word differences $d, d'$ and $x, x' \in X \cup \{e\}$. The target state of $(w, v)$ corresponds to the word difference $w^{-1}v$, and hence $D_K$ can be used to recognise equality between pairs of words which synchronously $K$-fellow travel.

Table 1: Properties of formal languages

| | Definition | True in |
|---|---|---|
| $P_\in$ | Membership of a language in $\mathcal{F}$ is decidable. Implied by $P_\emptyset \wedge P_{\cap \mathcal{R}}$. | $\mathcal{R}$, $\mathcal{C}$, $\mathcal{DC}$ $\mathcal{I}$, $\mathcal{R}ec$. |
| $P_\emptyset$ | Emptiness of a language in $\mathcal{F}$ is decidable. | $\mathcal{R}$, $\mathcal{DC}$, $\mathcal{C}$, $\mathcal{I}$. |
| $P_\infty$ | Finiteness/infiniteness of a language in $\mathcal{F}$ is decidable | $\mathcal{R}$, $\mathcal{C}$, $\mathcal{I}$. |
| $P_=$ | Equality of two languages in $\mathcal{F}$ is decidable | $\mathcal{R}$, not $\mathcal{C}$,$\mathcal{I}$. |
| $P_\cap$ | The intersection of two languages in $\mathcal{F}$ is in $\mathcal{F}$. | $\mathcal{R}$, $\mathcal{R}ec$, $\mathcal{R}en$, not $\mathcal{C}$, $\mathcal{I}$. |
| $P_\cup$ | The union of two languages in $\mathcal{F}$ is in $\mathcal{F}$. | Any AFL. |
| $P_\circ$ | The concatenation of two languages in $\mathcal{F}$ is in $\mathcal{F}$. | Any AFL. |
| $P_+$ | The positive closure $\cup_{n>0} L^n$ of a language $L$ in $\mathcal{F}$ is in $\mathcal{F}$. | Any AFL. |
| $P_*$ | The Kleene closure $\cup_{n\geq 0} L^n$ of a language $L$ in $\mathcal{F}$ is in $\mathcal{F}$. | Any AFL. |
| $P_c$ | For any $L \subset X^* \in \mathcal{F}$, $X^* \setminus L$ is in $\mathcal{F}$. | $\mathcal{R}$, $\mathcal{DC}$, $\mathcal{R}ec$, not $\mathcal{C}$, $\mathcal{I}$, $\mathcal{R}en$. |
| $P_{\cap \mathcal{R}}$ | The intersection of a language in $\mathcal{F}$ with a regular set is in $\mathcal{F}$. | Any trio. |
| $P_h$ | The image of a language in $\mathcal{F}$ under a homomorphism (the restriction to the language of a semigroup homomorphism mapping $\epsilon$ to $\epsilon$) is in $\mathcal{F}$. | Any full trio. |
| $P_{\epsilon-free\,h}$ | The image of a language in $\mathcal{F}$ under a homomorphism $h$ with $h(a) \neq \epsilon$, $\forall a \neq \epsilon$ is in $\mathcal{F}$. | Any trio. |
| $P_{\epsilon-bdd\,h}$ | (limited deletion) The image of a language in $\mathcal{F}$ under a homomorphism $h$ for which $h$ maps at most $K$ consecutive symbols in any string of $L$ to $\epsilon$ (for some constant $K$) is in $\mathcal{F}$. | Any trio. |
| $P_{h^{-1}}$ | If $h$ is a homomorphism from $X^*$ to $Y^*$, and $L \in \mathcal{F}$ is a language over $Y$, then $h^{-1}(L) \in \mathcal{F}$. | Any trio. |
| $P_{GSM}$ | The image of a language in $\mathcal{F}$ under the operation of a generalised sequential machine (GSM) mapping is in $\mathcal{F}$. A generalised sequential machine is defined as in [18] to be a (possibly non-deterministic) finite state automaton with an attached output string. Corresponding to each transition of the underlying automaton is a (possibly trivial) word which is appended to the output string; when a word is accepted, the word currently on the output string is output. Implies in particular $P_h$, $P_{\cap \mathcal{R}}$. | Any full trio. |
| $P_{\epsilon-free\,GSM}$ | The image of a language in $\mathcal{F}$ under the operation of an $\epsilon$-free GSM mapping (for which there is non-trivial output on every transition) is in $\mathcal{F}$. Implies in particular $P_{\epsilon-free\,h}$, $P_{\cap \mathcal{R}}$. | Any trio. |

A difference machine $D_K$ can also be constructed from the word differences of an asynchronous combing. However the set of fellow travelling pairs of words is not recognisable simply as the language of $D_K$ as a finite state automaton; the pairs of words need to be fed asynchronously through the machine. Nonetheless, with appropriate usage of this machine, it is decidable whether or not two words $w$ and $v$ asynchronously $K$-fellow travel. Further, for any given word $w$, the set of words $v$ which asynchronously $K$-fellow travel with $w$ is regular, as is the subset of all such words $v$ for which $w^{-1}v$ is equal to any particular word difference.

Suppose that $G$ has a synchronous combing $L$, with associated fellow traveller constant $K$. Let $D_K$ be the appropriate difference machine, and $\preceq$ the shortlex word order (for which $w \preceq v$ if $w$ is shorter than $v$, or the same length as $v$ and preceding it lexicographically). The subset $R$ of $L(D_K)$ consisting of pairs $(w, v)$ with $w \preceq v$ is then regular. If $L$ is regular, then so is

$$L_0 = \{v \in L : \nexists w \in L, ((w, v) \in R)\}$$

In this case, $L_0$ is a bijective, regular, synchronous combing for $G$. Hence every automatic group has a bijective, regular, synchronous combing.

The theoretical solubility of the word problem for an asynchronously combable group is, as we mentioned in section 2, an immediate consequence of Gersten's result that such a group has exponential isoperimetric inequality. When the combing $L$ is recursively enumerable, just as in [10] for automatic groups, a practical solution to the word problem can be described, which is based on the Turing machine which recognises $L$, and the associated difference machine $D_K$.

This algorithm to test that a given word $v$ in the generators is trivial in the group $G$ breaks down into three steps. The first step finds $w_e \in L$, with $w_e =_G e$, the second finds $w_v \in L$, with $w_v =_G v$, and the final step tests to see if $w_v =_G w_e$. Given any $w_0 = x_1 x_2 \ldots x_r \in L$, $w_e$ can be found as the last word in a sequence of words $w_0, w_1, \ldots w_r$ in $L$, such that $w_{i+1} x_{r-i} =_G w_i$, and, where $v = y_1 y_2 \ldots y_s$, $w_v$ can be found as the last word in a sequence $w_r, w_{r+1}, \ldots, w_{r+s}$, such that $w_{r+i} =_G w_{r+i-1} y_i$; thus the test is dependent on being able, given $u \in L$ and $x \in X$, to find $u' \in L$, with $u' =_G ux$. The set of all such $u'$ is found as the intersection of $L$ with the regular set of elements $u'$ which both $K$-fellow travel with $u$ and satisfy $u' =_G ux$.

For automatic groups this algorithm is easily seen to be quadratic (in the size of the generating set). But the proof of this fact depends on the fact that the lengths of the words $w_i$ can be controlled. This is implied by the existence of 'multipliers' (which is in turn implied for automatic groups by the closure of the family of regular sets under intersection) and the pumping lemma for regular sets; hence in general we cannot expect such a bound on the complexity.

## 5 Finite variations and change of generators

In this section we examine how moving to a group related by some kind of finite variation to a given group $G$ or changing the generating set for $G$ affects combability.

The following technical lemma will be useful.

**Lemma 5.1** *Suppose that $G$ is a group with finite generating sets $X$ and $Y$, and that $L$ is a language for $G$ over $X$. Suppose that $\phi$ is a map from $L$ into $Y^*$, with the property that, for constants $\lambda, C > 0$, for any $v \in L$,*

$$d_{G,Y}(\phi(v)(t), v(\lambda t)) < C, \quad \forall t.$$

*Then any fellow traveller properties satisfied by $L$ over $X$ are satisfied by $\phi(L)$ over $Y$.*

Note: it is not required that $\phi(v) =_G v$.

PROOF: We start by fixing a set of words $\{w_x : x \in X\}$ over $Y$, such that $w_x =_G x$, and a set of words $\{w_y : y \in Y\}$ over $X$, such that $w_y =_G y$. Let $M$ be the maximum length of any such $w_x$ or $w_y$.

Now suppose that $v, w$ are two words in $L$, and $\phi(v), \phi(w)$ their images in $\phi(L)$.

Suppose first that $v$ and $w$ asynchronously fellow travel over $X$, with relative speed function $h$, and fellow traveller constant $K$. Define $H : \mathbf{R} \to \mathbf{R}$ by $\lambda H(t) = h(\lambda t)$. Then

$$\begin{aligned}
d_{G,Y}(\phi(v)(t), \phi(w)(H(t))) &= d_{G,Y}(\phi(v)(t), v(\lambda t)) + d_{G,Y}(\phi(w)(H(t)), w(\lambda H(t))) \\
&+ d_{G,Y}(v(\lambda t), w(\lambda H(t))) \\
&\leq 2C + M d_{G,X}(v(\lambda t), w(h(\lambda t))) \leq 2C + MK.
\end{aligned}$$

So $\phi(v)$ and $\phi(w)$ asynchronously fellow travel over $Y$ with relative speed function $H$ and fellow traveller constant $2C + MK$.

Since $h' = H'$, we see that $\phi(v)$ and $\phi(w)$ boundedly asynchronously fellow travel over $Y$ if and only if $v$ and $w$ do so over $X$. If $h(t) = t$ for all $t$, then $H(t) = t$ for all $t$, and so further, when $v$ and $w$ synchronously fellow travel, so do $\phi(v)$ and $\phi(w)$.

We verify separately that $\phi(L)$ inherits right and left fellow traveller properties from $L$.

To see that $\phi(L)$ inherits right fellow traveller conditions from $L$, suppose that $v, w \in L$ and that $\phi(w) =_G \phi(v)$ or $\phi(w) =_G \phi(v)y$, for some $y \in Y$. Then $\phi(v) =_G vg$, for some $g \in G$ which has length at most $C$ as a word over $Y$, and hence length at most $MC$ as a word over $X$; similarly, for some $g'$, $\phi(w) =_G wg'$. So $d_{G,X}(v,w) < 2MC+1$, and so if $L$ satisfies a right fellow traveller condition with fellow traveller constant $K$, $v$ and $w$ must fellow travel in the appropriate way at distance at most $K(2MC+1)$. The above then implies that $\phi(v)$ and $\phi(w)$ similarly fellow travel (with an appropriate constant).

Now suppose that $L$ satisfies both right and left fellow traveller conditions. We want to show that $\phi(L)$ then satisfies left fellow traveller conditions. Suppose that $v, w \in L$ and that $\phi(w) =_G y\phi(v)$ for some $y \in Y$. Then $\phi(v) =_G vg$ and $\phi(w) =_G wg'$, for $g, g'$ as above, so $yvg =_G wg'$. Then the right and left fellow traveller properties of $L$ together imply that $v$ and $w$ fellow travel in the appropriate way at a distance of at most $K(2MC+1)$. Hence $\phi(v)$ and $\phi(w)$ similarly fellow travel. □

The results in the following proposition are mostly generalisations of results for automatic groups in [10]. Most of these results are also proved in a more general language theoretic setup (but not as general as ours) in [7] for asynchronous (but not synchronous) combings.

**Proposition 5.2 (Finite variations and change of generators)** *Let $\mathcal{F}$ be a family of languages, and $\tau$ one of the 6 classes of combings of section 2.*

**(a)** *Let $N$ be a finite normal subgroup of $G$. If $G$ has an $\mathcal{F}$-combing of type $\tau$, then so does $G/N$. Conversely, if $\mathcal{F}$ is closed under concatenation with a finite set, and $G/N$ has an $\mathcal{F}$-combing of type $\tau$, then so does $G$.*

**(b)** *Suppose that $\mathcal{F}$ is closed under concatenation with a finite set. Then if $G$ has an aynchronous, boundedly asynchronous or synchronous $\mathcal{F}$-combing so does any group $J$ of which $G$ is a subgroup of finite index.*

**(c)** *Suppose that $\mathcal{F}$ is closed under $\epsilon$-free GSM mappings, and limited deletion. Then if $G$ has an $\mathcal{F}$-combing $L$ of type $\tau$*

  **(i)** *$G$ has an $\mathcal{F}$-combing of type $\tau$ over any finite generating set,*

  **(ii)** *any subgroup $H$ of finite index in $G$ has an $\mathcal{F}$-combing of type $\tau$. If $L$ is bijective then so is the combing for $H$.*

Note that all the language theoretic conditions of this result are satisfied by $\mathcal{R}$, $\mathcal{C}$, $\mathcal{I}$, $\mathcal{R}ec$, $\mathcal{R}en$, and of course the family $\mathcal{U}$ of all languages with finite alphabet. Note also that part (b) of the result does not say anything about the bicombability of $J$.

PROOF:

**(a)** Let $\nu$ represent the natural homomorphism from $G$ to its quotient $G/N$.

Let $X = \{x_1, \ldots, x_n\}$ be a generating set for $G$, and $L$ an $\mathcal{F}$-combing of type $\tau$ for $G$ over $X$. Then the set $X' = \{x_1^\nu, \ldots x_n^\nu\}$ (in which no identifications are made between images of distinct elements of $X$ which are equal in the quotient group) generates $G/N$. Define $L'$ to be the set of images of elements of $L$ under $\nu$.

As a formal language $L'$ is equivalent to $L$. It remains to verify the appropriate fellow traveller conditions. Hence suppose that $v^\nu, w^\nu$ are elements of $L'$, where $v, w \in L$. To see that $L'$ inherits the right fellow traveller properties of $L$, note that if $w^\nu =_{G/N} v^\nu$, then $w =_G vn_1$ for some element $n_1 \in N$, and if $w^\nu =_{G/N} v^\nu x_i^\nu$, then $w =_G vx_i n_2$ for some element $n_2 \in N$. Since $N$ is finite, both $n_1$ and $x_i n_2$ can be expressed as words over $X$ of bounded length at most $|N|$. The fellow travel condition on $L$ then implies that $v$ and $w$ fellow travel in the appropriate way at distance at most $|N|K$. To see that $L'$ inherits the left fellow traveller properties of $L$, note that if $w^\nu =_{G/N} x_j^\nu v^\nu$, then $w =_G n_3 x_j v$, for some $n_3 \in N$, where $n_3 x_j$ can be written as a word of length at most $|N|$ over $X$. Hence $L'$ is a combing of the same type as $L$.

Now to prove the converse statement, suppose that $Y$ is a generating set for $G/N$, and that now $L'$ is an $\mathcal{F}$-combing of type $\tau$ for $G/N$. Let $Z$ be a set of elements of $G$ which maps bijectively to $Y$ under $\nu$. Then $Z \cup N$ generates $G$. Choose $L''$ to be the natural set of words over $Z$ which maps bijectively under $\nu$ to $L'$. $L''$ is equivalent to $L'$ as a formal language. It is clear that the concatenation of $L''$ with $N$ is a language in $\mathcal{F}$ for $G$ over $Z \cup N$.

Now suppose that $v = un$ and $w = u'n'$ are two words in the language $L''N$, with $u, u' \in L''$, $n, n' \in N$. If $w =_G v$, or $w =_G vn''$, for some $n'' \in N$, then $u^\nu =_{G/N} u'^\nu$. If $w =_G vz$, for some $z \in Z$, then $u^\nu y =_{G/N} u'^\nu$, where $y = z^\nu \in Y$. Hence right fellow traveller properties of $L'$ are inherited by $L''N$.

Similarly if $w =_J n''v$, then $u^\nu =_G u'^\nu$, and if $w =_J zv$, then $y^\nu u^\nu =_G u'^\nu$. Hence left fellow traveller properties of $L'$ are inherited by $L''N$.

**(b)** Let $L$ be an $\mathcal{F}$-combing for $G$ over an alphabet $X$.

Let $T$ be a set of left coset representatives for $G$ in $J$. Certainly the concatenation $L'$ of $L$ with $T$ is a language in $\mathcal{F}$ for $J$ over $X \cup T$. Now suppose that $v = ut$ and $w = u't'$ are two words in that language, with $u, u' \in L$, $t, t' \in T$. If $w =_J v$, then $u =_G u'$. If $w =_J vy$, for some $y \in X \cup T$, then $uu_{t,y} =_G u'$, where $u_{t,y} \in L$ is a word representing a Schreier generator of $G$ defined by the rule $ty =_G u_{t,y} t_{t,y}$, for $t_{t,y} \in T$ and hence (since the set of such generators is finite) has bounded length as a word over $X$. Hence right fellow traveller properties of $L$ are inherited by $L'$.

However left multiplication by elements of $T$ is not controlled, and so left traveller properties of $L$ are not clearly seen to be inherited by $L'$.

**(c) (i)** An asynchronous version of this result is proved in [7]; achieving a synchronous or boundedly asynchronous fellow traveller condition is harder. The argument is basically that of [10], but stated in more general language theoretic terms.

First suppose that $X$ is a generating set which does not contain the identity element. Then $G$ has an $\mathcal{F}$-combing over $X$ if and only $G$ has an $\mathcal{F}$-combing over $X \cup \{e\}$. One implication is obvious; a combing over $X$ is also a combing over $X \cup \{e\}$. On the other hand, given a combing $L$ in $\mathcal{F}$ over $X \cup \{e\}$, we construct (as in [10]) a language $L'$ over $X$ from $L$, by replacing (for some $m > 0$) every $m$-th occurrence of $e$ in a word by a particular length $m$ word $w_e$, and deleting all other occurrences of $e$. $L'$ is easily seen to be the image of $L$ under a GSM mapping. Although this GSM mapping is not $\epsilon$-free, we can find it as the composite of an $\epsilon$-free GSM mapping (which replaces all but the $m$-th occurrences of $e$ by some dummy symbol), and a homomorphism with limited deletion (of at most $m-1$ consecutive symbols). Hence $L'$ is in $\mathcal{F}$. The conditions of lemma 5.1 are satisfied with $\lambda = 1$ and $C = m$; hence $L'$ inherits the fellow traveller conditions of $L$.

Now let $X$, $Y$ be finite generating sets, and assume that $G$ has an $\mathcal{F}$-combing of type $\tau$ over $X$. By the above, we may without loss of generality assume that $e \in Y$. For each $x \in X$, we can find some word $w_x \in Y^*$ which is equal to $x$ in $G$. By appending copies of $e$ where necessary, we can ensure that all the words $w_x$ have the same length $m$. A language $L''$ over $Y$ is now defined as the image of $L$ under the homomorphism which maps $x$ to $w_x$;the condition $P_{\epsilon-freeGSM}$ implies $P_{\epsilon-freehom}$, and so ensures that $L''$ is in $\mathcal{F}$. The homomorphism expands the lengths both of words in $L$ and their subwords in $L$ by a constant factor $m$; hence lemma 5.1 applies.

If only asynchronous fellow traveller properties are required, it is not essential for $\mathcal{F}$ to satisfy $P_{\epsilon-bdd}$, since padding with copies of the symbol $e$ is no longer necessary; hence in this case it is sufficient for $\mathcal{F}$ to satisfy $P_{\epsilon-freehom}$.

**(ii)** Now suppose that $G$ has an $\mathcal{F}$-combing $L$ over a generating set $X$. Suppose that $T$ is a right transversal for $H$ in $G$, containing the identity. As above, the Schreier generators for $H$ are defined to be all the products $y_{tx} = txt_{tx}^{-1}$, for $x \in X$, and $t, t_{tx} \in T$, where $t_{t,x}$ is the representative in $T$ of the coset containing $tx$. An $\epsilon$-free generalised sequential machine $M$ can be defined, which first rewrites any word $x_1 \ldots x_k$ in $G$ as a product $y_1 \ldots y_k t$, for Schreier generators $y_i$ and $t \in T$, by reading from left to right and applying a succession of rewrite rules of the form $tx \to y_{tx}t_{tx}$, and then accepts $x_1 \ldots x_k$ and outputs $y_1 \ldots y_k$ in the case where $t =_G e$, but otherwise rejects $x_1 \ldots x_k$. Let $\phi$ be the associated GSM mapping. Then where $v = x_1 \ldots x_k$ and $\phi(v) = y_1 \ldots y_k$, we see that $d_{G,Y \cup T}(\phi(v)(i), v(i)) \leq 1$. Hence we can apply lemma 5.1 to see that $\phi(L)$ satisfies the same fellow traveller conditions over $Y \cup T$ that $L$ satisfies over $X$. Since any word in $\phi(L)$ is in fact a word over $Y$, and word differences between such words are therefore elements of $H$, a little thought shows that $\phi(L)$ satisfies the same fellow traveller conditions over the subset $Y$ of $Y \cup T$. It is straightforward to see that if $L$ is bijective then so is $\phi(L)$.

□

## 6 Products and extensions of combable groups

In this section we look at free and direct products, central extensions and (nearly) split extensions of $\mathcal{F}$-combable groups, for $\mathcal{F}$ a class of formal languages. As in the previous section, all the language theoretic conditions of these results are satisfied by $\mathcal{R}$, $\mathcal{C}$, $\mathcal{I}$, $\mathcal{R}ec$, $\mathcal{R}en$, and of course the family $\mathcal{U}$ of all languages with finite alphabets.

The situation for free and direct products is quite straightforward. For central extensions, a result of Neumann and Reeves ([20]) for biautomatic groups can be adapted for other language classes. Finally, a result of Bridson ([6]) deals with split extensions of combable groups over groups, such as abelian or word hyperbolic groups, which possess rather stable $\mathcal{F}$-combings; Bridson's result is already rather general, but it seems worthwhile, since it is clearly a very useful result, to state and prove it here in the slightly more general form which is easily obtained.

**Proposition 6.1 (Free products)** *Let $\mathcal{F}$ be a family of languages, and $\tau$ one of the six classes of combings described in section 2.*

*Suppose that $\mathcal{F}$ is closed under concatenation, Kleene closure, and intersection with regular sets, and that $G_1$ and $G_2$ have $\mathcal{F}$-combings of type $\tau$, Then $G_1 * G_2$ has an $\mathcal{F}$-combing of type $\tau$, provided that if $\tau$ consists of synchronous combings, the combings for $G_1$ and $G_2$ are bijective.*

PROOF: Where $L_1$ and $L_2$ are the languages of $G_1$ and $G_2$, over generator sets $X_1$ and $X_2$, and $L_1', L_2'$ are the subsets of those which contain no non-trivial representatives of the identity, then $L = (L_1' L_2')^*$ is an appropriate language for $G_1 * G_2$ over $X_1 \cup X_2$.

The fellow traveller properties of $L_1$ and $L_2$ ensure that each contains at most finitely many non-trivial representatives of $e$, and hence $L_1'$ and $L_2'$ can be found as the intersections of $L_1$ and $L_2$ with regular sets, and are in $\mathcal{F}$. Hence $L \in \mathcal{F}$.

Aysnchronous fellow traveller properties are easily seen to be inherited by $L$ from $L_1$ and $L_2$; for if two words $w, w'$ in $(L_1 L_2)^*$ are related by an equation of the form $w =_G w'$, $w =_G w'x$ or $w =_G xw'$, for $x \in X_1 \cup X_2$, then, for some $k$, $w = v_1 v_2 \ldots v_k$ and $w' = v'_1 v'_2 \ldots v'_k$, where $v_i =_G v'_i, v'_i x$ or $xv'_i$, and, for all but one $i$, $v_i =_G v'_i$. In the synchronous case, bijectivity ensures that for all but one $i$, $v_i = v'_i$.

Essentially this proof is given in [3] for automatic groups; there bijectivity can be assumed. □

**Proposition 6.2 (Direct products)** *Suppose that $\mathcal{F}$ is closed under concatenation, and that $G_1$ and $G_2$ have $\mathcal{F}$-combings $L_1$ and $L_2$. Then $L_1 L_2$ is an asynchronous $\mathcal{F}$-combing for $G_1 \times G_2$, and an asynchronous bicombing if $L_1$ and $L_2$ are asynchronous bicombings. If $L_1$ and $L_2$ are bijective, then so is $L_1 L_2$.*

*Where $L_1$ and $L_2$ are both synchronous or boundedly asynchronous then so is $L_1 L_2$ provided that $L_1$ is a near geodesic language.*

PROOF: The proof follows the argument of [10] for automatic groups, and uses the obvious decomposition of the Cayley graph for $G_1 \times G_2$ (with respect to a sensible generating set) into Cayley graphs for $G_1$ and $G_2$. The near geodesic condition ensures that, for $v, v' \in L_1, w, w' \in L_2$ with $vw =_G v'w'$, $vwx =_G v'w'$ or $yvw =_G v'w'$ (for generators $x, y$), the lengths of $v$ and $v'$ differ by a finite amount; for automatic groups, a language can always be found in which this condition holds. □

In fact this direct product result is a corollary of the asynchronous split extension result which is to come.

Neumann and Reeves' result on (virtually) central extensions of biautomatic groups of [20, 21] can be generalised to central extensions of $\mathcal{F}$-combable groups of various types without much difficulty. However it is not clear that the restrictive conditions will often be satisfied, for non-regular combings. Before we can state that result, we need some notation.

Let $H$ be a group with a combing $L$ and generating set $X$, and $A$ a finitely generated abelian group, on which $H$ acts as a group of automorphisms. An extension $G$ of $A$ by $H$ is then specified by a section $s : H \to G$ for $H$ in $G$, with $s(e) = e$, together with a cocycle $\sigma : H \times H \to A$, which satisfies the rule

$$s(h)s(h') =_G s(hh')\sigma(h, h')$$

Group elements can be written as products of the form $s(h)a$, and the above rule, together with the action of $H$ on $A$, defines multiplication between those products. The extension of $A$ by $H$ is *central* if $H$ acts trivially on $A$.

$\sigma$ is defined (as in [20]) to be an *L-regular cocycle* if
(a) the sets $\sigma(X, H)$ and $\sigma(H, X)$ are finite, and
(b) for each $x \in X$ and $a \in A$, the set $\{v \in L : \sigma(v, x) = a\}$ is regular.

**Proposition 6.3 (Central extensions)** , *a generalisation of [20], theorem A.*

*Let $H$ be a group with an $\mathcal{F}$-combing $L$ over a generating set $X$. Let $G$ be a central extension of $A$ by $H$, defined by a section $s : H \to G$ and a cocycle $\sigma : H \times H \to A$. Suppose that $\mathcal{F}$ is closed under GSM mappings and concatenation with a regular set. Suppose also that $\sigma$ is $L$-regular. Then $G$ has an asynchronous $\mathcal{F}$-combing $L''$. If $L$ is a bicombing, then so is $L''$.*

*If $L$ is synchronous or boundedly asynchronous, then so is $L'$, provided that $L$ is a near geodesic language.*

PROOF: The proof is basically that of Neumann and Reeves, translated into the appropriate language.

Define $Y$ to be the set of all elements of $G$ of the form $y_{x,h} = s(x)\sigma(h,x)^{-1}$ for $h \in H, x \in X$. The finiteness of $\sigma(X,H)$ ensures that $Y$ is finite. Let $Z$ be a generating set for $A$, and $L_A$ a regular synchronous bicombing for $A$.

For each $w = x_1 x_2 \ldots x_n \in L$,

$$s(w) =_G s(x_1)s(x_2)\sigma(x_1,x_2)^{-1} s(x_3)\sigma(x_1x_2,x_3)^{-1} \ldots s(x_n)\sigma(x_1x_2\ldots x_{n-1},x_n)^{-1}$$

Hence $s(w)$ is represented by the word

$$w_0 = y_{x_1,e} y_{x_2,x_1} y_{x_3,x_1x_2} \cdots y_{x_n,x_1x_2\ldots x_{n-1}}$$

Let $L'$ be the set of all words of that form, and define $L''$ to be the language $L'L_A$. Then $L''$ is clearly a language for $G$.

The fellow traveller properties of $L''$ can be seen to be inherited from $L$ using arguments similar to those of [20], as follows.

First suppose that $w, w' \in L$ and let $w_0, w'_0 \in L'$ represent $s(w), s(w')$. Then for any $t$, $w_0(t)$, and $w'_0(t)$ are elements of $L'$ representing $s(w(t))$ and $s(w'(t))$. Hence if $w_0 =_G w'_0$ or $w_0 =_G w'_0 y$ for some $y \in Y$, right fellow traveller properties of $L$ imply that $w_0$ and $w'_0$ fellow travel in the appropriate way. Left fellow traveller properties of $L$ (if they exist) imply that $w_0$ and $yw'_0$ fellow travel when $w_0 =_G yw'_0$

For $v, v' \in L_A$, fellow traveller properties between pairs of words $w_0 v$ and $w'_0 v'$ which satisfy $w_0 v =_G w'_0 v'$, $w_0 v =_G w'_0 v' t$ or $w_0 v =_G t w'_0 v'$ (for $t \in Y \cup Z$) are now deduced, as in [20], from the fellow traveller properties of $L''$ and $L_A$ separately. A condition that $L$ is a near geodesic language ensures that the difference in lengths of $v_0$ and $v'_0$ remain bounded; for the automatic groups considered by [20] this property can always be assumed to hold.

It remains to show that $L''$ is in $\mathcal{F}$. We show this by effectively translating the construction of [20] into language theoretic terms. Since $\mathcal{F}$ is closed under concatenation with a regular set, it is enough to show that $L'$ is in $\mathcal{F}$. We construct a GSM mapping which maps $L$ to $L'$. For each $x \in X$ and each $a \in \sigma(G, X)$, let $W_{x,a}$ be the finite state automaton accepting $\{w \in L : \sigma(w, x) = a\}$, and let $S_{x,a}$ be its set of states; there are finitely many such sets. We form a (non-deterministic) generalised sequential machine $M$ with input alphabet $X$ and output alphabet $Y$ as follows. The state set of $M$ consists of the elements of the Cartesian product $\prod S_{x,a}$, together with a failure state $F$. For any input symbol $x$, any $a \in \sigma(G, X)$, and any state $s$ for which the $x, a$ component is an accept state, there is an arrow in $M$ taking $s$ to the state of which each component is the target under $x$ of the corresponding component of $s$, and such that the symbol $y_{x,a}$ is appended to the output. $M$ maps $L$ to $L'$. □

The next result is a (slight) generalisation of a result from [6], where it is the 'main lemma'; Bridson's result constructs a bijective asynchronous combing for a split extension of two groups out of bijective asynchronous combings for the two groups; in fact the construction also works when the extension is nearly split, and bijectivity is not strictly necessary. This seems very much to be an asynchronous result.

We need a little notation before we can state the result.

Following [22], where $G$ is a group with normal subgroup $N$, we say that $G$ is an $(r, s)$-split extension of $N$ by $G/N$ if for some $H \subseteq G$, $NH$ has index at most $r$ in $G$, and $|N \cap H|$ divides $s$. If $G$ is an $(r, s)$ extension of $N$ for some $r, s$ we say that $G$ is a *nearly split* extension; if $r = s = 1$, then $G$ is a split extension.

Suppose that groups $N$ and $H$ have asynchronous combings $L_N, L_H$ over generating sets $X, Y$, and suppose that $H$ acts as a group of automorphisms of $N$; for $n \in N$, denote by $n^h$ the image of $n$ under $h$. For $n \in N$, let $v_n = x_1 x_2 \ldots x_m$ be a word in $L_N$ which represents $n$. We call a word $w$ an *image* of $v_n$ under $h$ if $w$ is of the form $v_{x_1^h} v_{x_2^h} \ldots v_{x_m^h}$ such that, for each $i$, $v_{x_i^h}$ is a word in $L_N$ representing $x_i^h$.

**Proposition 6.4 ((Nearly) split extensions)** , *a slight generalisation of Bridson's lemma, [6].*

*Let $G$ be a nearly split extension of groups $N \triangleleft G$ and $H$. Suppose that $N$ and $H$ have asynchronous $\mathcal{F}$-combings $L_N$ and $L_H$ over generating sets $X$ and $Y$, for some language family $\mathcal{F}$ which is closed under concatenation. If $H \cap N$ is non-trivial, suppose further that $L_N$ is an asynchronous bicombing.*

*Suppose also that for all $n \in N$, all $y \in Y$, and all $v_n$ and $v_{n^y}$ in $L$ representing $n$ and $n^y$ respectively, $v_{n^y}$ asynchronously fellow travels with each of the images of $v_n$ under $y$ (Condition (∗)).*

*Then $G$ has an asynchronous combing, $L_G$. If $L_H$ and $L_N$ are bijective and the extension is split, then $L_G$ is also bijective.*

Note that it is enough to require that for each $n, h$, some $v_{n^h}$ asynchronously fellow travels with some image of $v_n$ under $h$; the rest of that condition follows from the fellow traveller properties of $L_N$.

PROOF: By proposition 5.2 (b), we lose no generality by assuming that $G$ is equal to the product $HN$. In this case, we define $L_G = L_H L_N$.

Throughout the proof, we use the notation $v_n, v'_n$ for elements of the language $L_N$ which represent $n$, and $w_h, w'_h$ for elements of $L_H$ representing $h$.

Every element of $G$ can be written in the form $hn$ for some $n \in N$, and some $h \in H$, so $L_H L_N$ is clearly a language for $G$. For $n_1, n_2 \in N$, $h_1, h_2 \in H$,

$$h_1 n_1 . h_2 n_2 =_G h_1 h_2 n_1^{h_2} n_2$$

Let $\Gamma_H, \Gamma_N, \Gamma_G$ be the Cayley graphs of the three groups. Of course $\Gamma_H$ and $\Gamma_N$ embed in $\Gamma_G$; in fact each embeds many times (for any vertex of $\Gamma_G$ an embedding may be chosen which places the basepoint of $\Gamma_N$ or $\Gamma_H$ at that vertex). If $u$ and $u'$ are paths of $\Gamma_G$ which lie in the same embedding of $\Gamma_H$ or $\Gamma_N$, then any fellow traveller condition satisfied by them in that graph clearly also holds between them in $\Gamma_G$.

Suppose that $n \in N$, $h \in H$, and let $x \in X \cup \{e\}$. Let $u_1 = w_h v_n$ and $u_2 = w'_h v'_{nx}$ be words in $L_G$, for which $w_h$ and $w'_h$ represent $h$ in $L_H$, and $v_n$ and $v_{nx}$ represent $n$ and $nx$ in $L_N$. The portions of $u_1$ and $u_2$ along $w_h$ and $w'_h$ asynchronously fellow travel within an embedded $\Gamma_H$, and the portions from $h$ to $hg$ and $hgx$ along $v_g$ and $v_{gx}$ asynchronously fellow travel within an embedded $\Gamma_G$ (based at $h$). Hence $u_1$ and $u_2$ asynchronously fellow travel.

Now choose $y \in Y$. Since $hny = hyn^y$, $hny$ is represented in $L_G$ by any word $u_3 = w_{hy} v_{n^y}$, for which $w_{hy} \in L_H$ represents $hy$ and $v_{n^y} \in L_N$ represents $n^y$. The portions of $u_1$ and $u_3$ along $w_h$ and $w_{hy}$ asynchronously fellow travel in $\Gamma_G$ because they do in $\Gamma_H$. Now the portion of $u_3$ from $hy$ to $hny$ along $v_{n^y}$ asynchronous fellow travels with any path $u'$ labelled by an image of $v_{n^y}$ under the action of $y^{-1}$ running from $h$ to $hn$; edges labelled $y$ run from any such path to the second part of $u_3$. The condition $(*)$ implies that any such path $u'$ asynchronously fellow travels (within an embedded copy of $\Gamma_N$) with the path from $h$ to $hn$ labelled by $v_n$, that is, with the second part of $u_1$; hence we have verified that $u_1$ and $u_3$ asynchronously fellow travel.

In the case where $N \cap H$ is trivial, any element of $G$ has a unique representation as a product of the form $hn$, for $h \in H$, $n \in N$, and so the above analysis covers all representatives in $L_G$ of elements $g, g'$ for which $g' =_G gx$ for $x \in X \cup \{e\}$ or $g' =_G gy$ for $y \in Y$. If $N \cap H$ is a non-trivial finite group $F$, the proof is complete once we have verified that, for $h, h' \in H$ and $n, n' \in N$, with $hn =_G h'n'$, words of the form $v_h w_n$ and $v_{h'} w_{n'}$ must asynchronously fellow travel. In this case the equation $hn =_G h'n'$ implies that $h' = hf$, $n' = f^{-1}n$ for some $f \in F$. Since $F$ is finite, $f$ has bounded length as a product of generators, and hence the fact that $v_h$ and $v_{h'}$ asynchronously fellow travel follows from the fellow traveller condition on $H$; similarly, $v_n$ and $v_{n'}$ asynchronously fellow travel, provided that $L_N$ is an asynchronous bicombing. □

To deduce a synchronous, rather than asynchronous, fellow traveller property, it seems that we need, in addition to synchronous fellow traveller properties in $L_N$ and $L_H$,
(a) that for all $y \in Y$, $x \in X$, $x^y$ is a generator
(b) a synchronous version of $(*)$.
Furthernore, to get conditions for any type of bicombing seems to be more complicated. For notice that, for $x \in X$, $xhg = hx^h g$. 'Left' fellow traveller conditions are more easily satisfied by the language $L_G L_H$, but this does not as easily have good 'right' fellow traveller properties.

Bridson observes in [6] that the condition $(*)$ is satisfied by any geodesic combing of a word hyperbolic group; these combings are regular. He also shows that any finitely generated abelian group has an asynchronous combing which satisfies $(*)$. For the torsion free abelian group $\mathbf{Z}^n$ this language is found by embedding the Cayley graph in the obvious way in $\mathbf{R}^n$, and selecting for each element $g$ the lexicographically earliest of the set of words which remain closest to the strange line path in $\mathbf{R}^m$ from $e$ to $g$; and in the case $n = 2$, this language is proved (in [7]) to be an indexed language. In both the above situations, the combings are in fact bicombings, and so the result can be applied to nearly split extensions.

## 7 Combings for nilpotent and soluble groups

In [7], Bridson's lemma was applied to show that split extensions of the form $\mathbf{Z}^2 \rtimes \mathbf{Z}$ have asynchronous indexed ($\mathcal{I}$) combings, and hence to prove that the same is true of $\pi_1(M)$ for any compact, geometrisable 3-manifold $M$. Of course the construction used by Bridson has much more general application. We shall briefly observe a few consequences in this section. We see that many nilpotent, and also non-nilpotent soluble groups can be expressed as (nearly) split extensions over $\mathbf{Z}^n$ and so have asynchronous combings. In some cases the combings are already proved (by the results of [7]) to be indexed languages. In general, it seems likely that the language for $\mathbf{Z}^n$ described in the last paragraph of section 6 is both an indexed language and a real time language (see, for example, [23]) for all values of $n$; in that case the combings described in this section would also lie in these language classes. These questions are explored more fully in [13].

Where $G$ is a non virtually abelian, nilpotent group, we cannot expect to do much better than find an indexed combing for $G$. For it is proved in [10] that no such group $G$ can have a regular combing, and in [8] (and [7]) that $G$ cannot have a bijective combing which is a bounded language. Using the latter result, it is proved in [7] that $G$ cannot have a bijective, context-free combing (here bijectivity is essential to the proof, which deduces polynomial growth of the combing from the polynomial growth of $G$).

The $n$-th Heisenberg group

$$G_n = \langle a_i, b_i (i = 1, \ldots n), c \mid c = [a_i, b_i], [a_i, c] = [b_i, c] = 1 \\ [a_i, a_j] = [a_i, b_j] = 1 (i \neq j) \rangle,$$

which is class 2 nilpotent, is a split extension of the free abelian group $\langle c, a_i b_i\ (i = 1, \ldots n) \rangle$ by the free abelian group $\langle a_1, \ldots a_n \rangle$. Hence this has a bijective, asynchronous combing.

Let $U_n$ be the group of $n \times n$ uni-uppertriangular integer matrices. This is nilpotent of both nilpotency and derived class $n - 1$; the terms of both the derived and lower central series are subgroups isomorphic to $U_{n-1}, U_{n-2}, \ldots, U_2 \cong \mathbf{Z}$. We can express $U_n$ as a split extension of $\mathbf{Z}^{n-1}$ (found as the matrices with 1's on the diagonal, and 0's elsewhere except in the right hand column) and $U_{n-1}$ (found as the matrices with 0's in the right hand column, except for the diagonal entry). Hence, by a clear induction argument $U_n$ has an asynchronous combing.

Let $N = N_{k,2}$ be the free nilpotent class 2 group on $k$ generators, with presentation

$$\langle x_1, \ldots x_k | [[x_i, x_j], x_k], \forall i, j, k \rangle$$

Then $N_{k,2}$ is a split extension of the normal abelian subgroup generated by $x_k$ and all commutators of the form $[x_k, x_j]$ by the subgroup generated by $x_1, \ldots x_{k-1}$. Hence (by induction) $N_{k,2}$ has an asynchronous combing. However, for $c > 2$, the free nilpotent class $c$ group $N_{k,c}$ on $k$ generators $x_1, \ldots, x_k$ does not split over an abelian subgroup; hence it is not clear whether or not $N_{k,c}$ has a combing for $c > 2$.

In fact, using Bridson's lemma, the examples for $N_{k,2}$ and the Heisenberg groups, and the finite variation results of proposition 5.2, many nilpotent groups can be shown to have asynchronous combings. The question is investigated further in [13]; in particular combings are found for all class 2 nilpotent groups with 2 or 3 generators or cyclic commutator subgroup.

We can use a result of Robinson ([22]) to show that many soluble groups which are far from being nilpotent are also in this class. Robinson's result is rather general; all we really need is the version given by Segal in [24], namely that if $A$ is a finitely generated, free abelian normal subgroup of a group $G$ such that $G/A$ is finitely generated and nilpotent, and such that $C_A(G) = 1$, then $G$ nearly splits over $A$; in fact a subgroup of finite index in $G$ is a split extension of $A$.

Using this result we can prove the following.

**Proposition 7.1** *If $G$ is polycylic, metabelian and torsion free with centre disjoint from $G'$, then $G$ has an asynchronous combing. When $G'$ has rank 2, the* *combing is an indexed language.*

PROOF: That $G$ is polycyclic and metabelian implies that $G'$ is finitely generated. By Robinson's theorem, $G$ nearly splits over $G'$, and we can apply Bridson's lemma. $\square$.

As an example, we have the group

$$\langle x, y, z \mid yz = zy, y^x = z, z^x = yz \rangle$$

which is certainly not automatic (it has exponential isoperimetric inequality, see [10], theorem 8.1.3).

# References

[1] Alfred V. Aho, Indexed grammars - an extension of context-free grammars, J. Assoc. Comp. Math (15) 1968, 647–671.

[2] Alfred V. Aho, Nested stack automata, J. Assoc. Comp. Mach. 16 (1969) 383–406.

[3] G. Baumslag, S.M. Gersten, M. Shapiro, H. Short, Automatic groups and amalgams, Journal of Pure and Applied Algebra 76 (1991) 229–316.

[4] G. Baumslag, D. Gildenhuys, R. Strebel, Algorithmically insoluble problems about finitely presented solvable groups, Lie and associative algebras. I, Journal of Pure and Applied Algebra 39 (1986) 53–94.

[5] Martin R. Bridson, On the geometry of normal forms in discrete groups, Proc. London Math. Soc. (3) 67 (1993) 596–616.

[6] Martin R. Bridson, Combings of semidirect products and 3-manifold groups, Geometric and Functional Analysis 3 (1993) 263–278.

[7] M.R. Bridson and R. H. Gilman, Formal language theory and the geometry of 3-manifolds, in Commentarii Math. Helv. 71 (1996) 525–555.

[8] M.R. Bridson and R. H. Gilman, On bounded languages and the geometry of nilpotent groups, Combinatorial and geometric group theory, LMS Lecture Notes Ser. 204, ed. Duncan, Gilbert and Howie, 1–15.

[9] J. W. Cannon,The combinatorial structure of cocompact discrete hyperbolic groups, Geom. Dedicata 16 (1984) 123–148.

[10] David B. A. Epstein J.W. Cannon, D. F. Holt, S. Levy, M. S. Patterson and W. Thurston,Word processing in groups, Jones and Bartlett, 1992.

[11] S. M. Gersten, Bounded cocycles and combings of groups, Internat. J. Algebra Comput. 2 (1992) 307–326.

[12] Robert H Gilman, Formal languages and infinite groups, DIMACS Ser. Discrete Math. Theoret. Comput. Sci., 25, AMS 1996, 27–51.

[13] Robert F. Gilman, Derek F. Holt and Sarah Rees, Combings for nilpotent groups, preprint in preparation.

[14] Seymour Ginsburg, The mathematical theory of context-free languages, McGraw-Hill, New York etc., 1966.

[15] S Ginsburg and S Greibach, Abstract families of languages, IEEE Conf. Record of Eighth Ann. Symp. on Switching and Automata Theory (1967), 128–139.

[16] Seymour Ginsburg and Edwin H. Spanier, Bounded algol-like languages, Trans. A.M.S. 113 (1964) 333–368.

[17] T. Hayashi, On derivation trees of indexed grammars, Publ. RIMS Kyota Univ. 9 (1973) 61–92.

[18] John E. Hopcroft and Jeffrey D. Ullman, Introduction to automata theory, languages and computation, Addison-Wesley, 1979.

[19] O.G. Kharlampovich, A finitely presented group with unsolvable word problem (Russian), Izv. Akad. Nauk SSSR Ser. Mat. 45 (1981) 852–873.

[20] Walter D. Neumann and Lawrence Reeves, Regular cocycles and biautomatic structures, Internat. J. Algebra Comput. 6 (1996) 313–324.

[21] Walter D. Neumann and Lawrence Reeves, Central extensions of word hyperbolic groups, Annals of Math. 145 (1997) 183–192.

[22] Derek J.S. Robinson, Splitting theorems for infinite groups, Symposia Mathematica 17, Convegni del Novembre e del Dicembre 1973, Academic Press, 1976, 441–470.

[23] Arnold L. Rosenberg, Real-time definable languages, J. Assoc. Comput. Mach. 14 (1967) 645–662.

[24] Daniel Segal, Polycyclic groups, Cambridge University Press, 1983.